\documentclass[a4paper,fleqn]{cas-sc}

\usepackage[numbers]{natbib}
\newtheorem{propo}{Proposition}[section]
\newtheorem{lemma}[propo]{Lemma}
\newtheorem{example}{Example}
\newtheorem{remark}{Remark}

\begin{document}

\let\WriteBookmarks\relax
\def\floatpagepagefraction{1}
\def\textpagefraction{.001}

\shorttitle{\'A. Carmona, A. M. Encinas, M.J. Jim\'{e}nez and \'A. Samperio}    

\shortauthors{\'A. Carmona et~al.}  


\title [mode = title]{Stable  recovery of  piecewise constant conductance on spider networks}

\tnotemark[] 

\tnotetext[1]{This work has been partially supported by the Spanish
Research Council (Ministerio de Ciencia e Innovaci\'on) under project 
PID2021-122501NB-I00 and by the Universitat Polit\`ecnica de Catalunya under funds AGRUPS-UPC 2022 and 2023. \'A. Samperio was supported by a FPI grant of the Research Project PGC2018-096446-BC21 (with the help of the FEDER Program).}

\author[1]{\'A. Carmona}[orcid=0000-0001-7713-1066]
\ead{angeles.carmona@upc.edu}
\affiliation[1]{organization={Universitat Polit\`ecnica de Catalunya - BarcelonaTech (UPC),  Department of Mathematics},
            city={Barcelona},
            country={Spain}}
\author[1]{A. M. Encinas}[type=editor, orcid=0000-0001-5588-0373]
\ead{andres.marcos.encinas@upc.edu}

\author[1]{M.J. Jim\'{e}nez}[type=editor, orcid=0000-0003-3502-462X]
\cormark[1]
\ead{maria.jose.jimenez@upc.edu}

\author[2]{\'A. Samperio}[type=editor, orcid=0000-0002-1889-8889]
\ead{alvaro.samperio@uva.es}

\affiliation[2]{organization={Departamento de \'Algebra, An\'alisis Matem\'atico, Geometr\'\i a y Topolog\'\i a and IMUVA,  Universidad de Valladolid},
            city={Valladolid},
            country={Spain}}

\cortext[1]{Corresponding author.\\E-mail addresses: }
 %---------------------------------------fin Datos personales------------------------------------------

\begin{abstract}
We address the discrete inverse conductance problem for well-connected spider networks; that is,  to recover the conductance function on a well-connected spider network from the Dirichlet-to-Neumann map. It is well-known that this inverse problem is exponentially ill-posed, requiring the implementation of a regularization strategy for numerical solutions. Our focus lies in exploring whether prior knowledge of the conductance being piecewise constant within a partition of the edge set comprising few subsets enables stable conductance recovery.  To achieve this, we propose formulating the problem as a polynomial optimization one, incorporating a regularization term that accounts for the piecewise constant hypothesis. We show several experimental examples in which the stable conductance recovery under the aforementioned hypothesis is feasible.

\end{abstract}


\maketitle

\section{Introduction and Preliminaries}

The inverse conductivity problem, known as   Calderon's inverse problem,
has a long history and numerous papers have been published in its continuous version. However, 
the discrete version of the problem, understood as its  approach in an electrical network, is more recent and has received comparatively less attention (see \cite{A88,AV05,ACE15-2,ACE15,BDG08,BDGM11,BDM10,BU97,C80,CIM98,R16,SU87}). 

The problem at hand involves an unknown conductivity that needs to be determined and possibly reconstructed using boundary measurements of current and voltage. This intriguing challenge has garnered significant attention due to its wide range of applications in diverse fields, including noninvasive medical imaging, which stands as one of the most complex and compelling areas of interest (see \cite{H05,CIN99,PHMM19,SYXRX21}).

Although the conductivity could be uniquely reconstructed from the data, we cannot assure stability as the problem is severely ill-posed. For this reason, various authors have reformulated the problem by considering the recovery of conductivity when some {\it a priori} information is known. For example, Alessandrini and Vessella proved in \cite{AV05}  that if it is {\it a priori} known that the conductivity is piecewise constant with a bounded number of unknown values, the problem becomes Lipschitz stable. Other authors have considered  regularization methods, mainly of Tikhonov type (see \cite{EHN96,LMP03,R16}). 

In the discrete setting; that is, in the case of electrical networks, Curtis and Morrow showed uniqueness of the inverse conductance problem for critical planar networks (see \cite{CM90,CM91,CM00,CMM94}). The problem has also been studied for trees (see \cite{GR22}) and for $n$-dimensional grids (see \cite{ACEM16}). Among the more used critical planar networks are the well-connected spider networks, which have been studied both as a discretization of a continuous media  to approximate the continuous problem and as the workspace to solve the discrete inverse problem (see \cite{ACE15-2,BDM10,CMM94}). Other authors have used grids to perform the discretization approach (see \cite{BDGM11,EdG11}).

Some of the authors obtained  an explicit formula for the recovery of the conductance on a well-connected spider network in \cite{ACE15-2}. In this work, we de\-mons\-trate that, in practice, the use of this formula leads to the inherent instabilities in the inverse problem, even for relatively small-sized networks. This poor performance explains why the use of these networks in medical applications is restricted to networks with fewer than 16 nodes on the boundary (see \cite{SYXRX21,Z20}).

Within this scenario, we investigate whether  imposing a hypothesis in the discrete case that mimics the piecewise constant conductivity assumption leads  to a stable recovery of the conductance. Under this hypothesis, we propose formulating the inverse problem as a polynomial optimization problem, where we penalize deviations from the assumption. This naturally leads us to  a regularization approach {\it \`a la Tikhonov}. 

Finally,  we present numerical experiments for the recovering of piecewise constant conductance on spider networks. This method allows us to obtain regularized solutions that are good approximations and support the stability of our approach. Moreover, our experiments  are also in agreement with  the results obtained in \cite{M01,R06} for the continuous case,  showing that the Lipschitz constant grows exponentially with the number of  subsets in which the conductance is constant.

Throughout the paper, we use the common definitions and notations on finite networks for the discrete inverse problem. In particular we mainly follow the terminology and techniques in \cite{ACE15} and also 
in \cite{ACE15-2,ACEM16}.

Consider a finite network denoted as $\Gamma$, where $\Gamma$ is defined as $(V, c)$. This network is essentially a  finite, connected graph comprising a set of vertices represented by $V$ and characterized by its conductance, denoted as $c$.

We use the notation ${\cal C}(F)$ to represent the set of functions defined on a subset $F \subseteq V$, which  can be  identified with $\mathbb{R}^{|F|}$, where $|F|$ stands for the cardinality of $F$. The boundary of the set $F$, denoted as $\delta(F)$, comprises the vertices in the complement of $F$, which is represented as $F^c$, and these vertices are connected to at least one vertex within the set $F$. Additionally, the closure of the set $F$ is defined as $\bar F = F \cup \delta(F)$.

When we assign labels to the vertices in the set, we represent the matrix related to the Laplacian operator as ${\sf L}.$
In addition, fixed the subset $F\subset V$, ${\sf N}$ denotes the  matrix associated with the Dirichlet-to-Neumann map, that in this setting in known as the {\it response matrix} corresponding to $F$. The characterization of those matrices of order $|\delta(F)|$ that can serve as response matrix of a network was carried out  in \cite{CIM98}. 

For  ${x,y\in V}$ (respectively $x,y\in \delta(F)$), $L(x,y)$ (respectively $N(x,y)$) is the entry of ${\sf L}$ (res\-pec\-ti\-ve\-ly ${\sf N}$) corresponding to  vertices $x$ and $y$. More generally, given  a pair of disjoint subsets $P,Q\subset V$, we consider the submatrix of  ${\sf L}$, ${\sf L}(P;Q)=\big(L(x,y)\big)_{(x,y)\in P\times Q}$. Similarly, in the case where we have subsets $P$ and $Q$ that are both within the boundary set $\delta(F)$, we consider a submatrix within the response matrix denoted as ${\sf N}(P;Q)$, which is formed by   $\big(N(x,y)\big)_{(x,y)\in P\times Q}$.

When the set of vertices $V$ is equal to the closure of $F$, and there are no adjacent boundary vertices, we refer to the network $\Gamma=(\bar F,c)$ as a {\it circular planar} network. This term is used when the underlying graph can be placed inside a closed disc, denoted as $D$, in the plane in such a way that the vertices within set $F$ are located inside the interior of $D$ ($\stackrel{\circ}{D}$), and the vertices in the boundary set $\delta(F)$ are positioned along the circumference of $D$ ($\partial D$). Additionally, the vertices in $\delta(F)$ can be assigned labels in a clockwise order.

In the case of a special type of circular planar graph known as the {\it critical} circular planar graph, you can determine the conductance uniquely by examining its response matrix, as demonstrated in \cite[Theorem 4]{CIM98}. As a result, it becomes meaningful to consider the inverse problem of recovering the conductance for these types of graphs. One noteworthy subset of critical circular planar networks is the well-connected spider networks, which were initially introduced in \cite{CM00} due to their exceptional characteristics.

\begin{figure}[htb!]
  \begin{center}
    \includegraphics[width=0.45\textwidth]{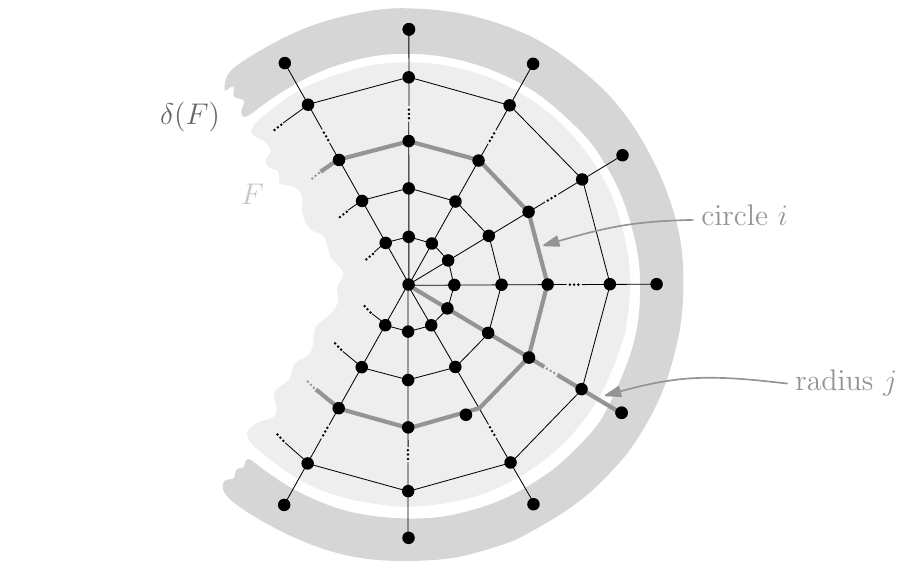}
  \caption{Representation of a spider network.} 
  \end{center}
  \label{fig:spidernetwork}
\end{figure}

A {\it well-connected spider network} with $\ell\ge 0$ circles and $m=4\ell+3$ radii can be described as a circular planar network, denoted as $\Gamma=(\bar F,c)$. In this network, there are a total of $m$ boundary vertices represented by $\delta(F)=\{1,\ldots,m\}$, and they are arranged in the circular order as defined by $\partial D$, the boundary of the circle $D$. The distribution of vertices within $F$ follows a specific pattern: Start by placing a vertex at the center of the circular boundary $\partial D.$
Draw straight lines, referred to as {\it radii}, from the central vertex to each of the boundary vertices.
Now, for the internal structure of the network: Create $\ell$ distinct concentric circumferences, all contained within the interior of $D$. Each of these concentric circles is referred to as a {\it circle}. 
Place a vertex at the intersection point of every circle and radius.
The network's edges are determined by these radii and circles, as shown in Figure \ref{fig:spidernetwork}.

In \cite{ACE15}, an explicit iterative process which recovers the conductance of a well-connected spider network from its Dirichlet-to-Neumann map was introduced. In the next section we show the features of this algorithm as well as the limitations.

\section{Ill-posedness: limitations of the explicit recovery process}\label{exact}

The aim of this section is to prove that, in the discrete setting, the inverse conductance recovery problem  is 
intrinsically ill-posed. As we  show in the following tests, for well-connected spider networks  if we compute the Dirichlet-to-Neumann map, and then we apply the algorithm in \cite{ACE15}, we recover a con\-duc\-tan\-ce that, when the number of  boundary vertices is high, widely differs  from the one of the original network. This is due to the ill-posedness of the problem: despite the algorithm is based on explicit formulas, any error in the entries of the Dirichlet-to-Neumann matrix (which are stored with finite precision) could be  amplified several orders of magnitude in the algorithm.

\begin{example}\label{exacto} For $m=7, 11, 15, 19, 23, 27, 31$ and $35$, we start from the well-connected spider network with $m$ radii and constant conductance $c=1$. In all cases, we compute the response matrix ${\sf N}$ of the network, and from it  we recover the conductance $c'$ using the explicit formulas from \cite{ACE15}. The algorithm has been implemented in Matlab. 
\end{example}

\begin{figure}[h!]
\begin{center}
\includegraphics[scale=0.6]{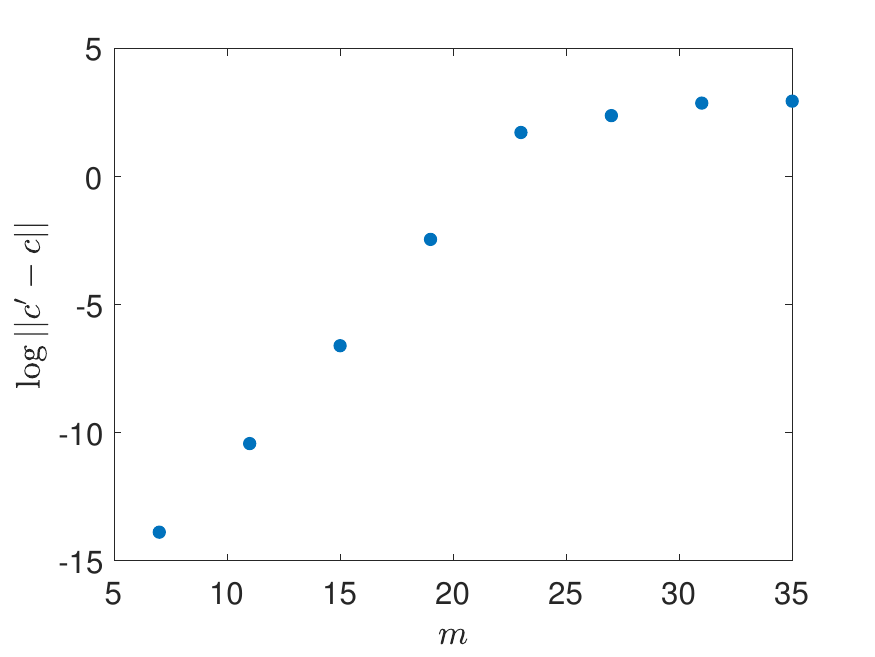}
\caption{Logarithm of the error in the recovered conductance.}
\label{exacto8}
\end{center}
\end{figure}

In Figure \ref{exacto8} we show the logarithm of the error in the recovered conductance in the Euclidean norm, $\log||c'-c||$, for all values of $m$, where $\log$ stands for the decimal logarithm. Moreover, Table \ref{tabla_exact} displays the error on the conductances. We see that the error is almost zero for $m=7$ and increases approximately exponentially with $m$ from $m=7$ to $m=23$. The error keeps increasing with $m$ for $m\geq 23$. 

\begin{table}[h!]
\caption{\label{tabla_exact}Error in the recovered conductance.}
{{\begin{tabular}{@{}lcccccccc}
\\[-3ex]
\noalign{\hrule height 1pt} 
$m$&7&11&15& 19&23&27&31 &35\\
\hline
$\log||c'-c||$&$1 \cdot 10^{-14}$&$4 \cdot 10^{-11}$&$3 \cdot 10^{-7}$&$4 \cdot 10^{-3}$&$5 \cdot 10^{2}$&$2 \cdot 10^{3}$&$7\cdot 10^{3}$&$9 \cdot 10^{3}$\\
\noalign{\hrule height 1pt} 
\end{tabular}}}
\end{table}

\begin{figure}[h!]
\begin{center}
\includegraphics[scale=0.45]{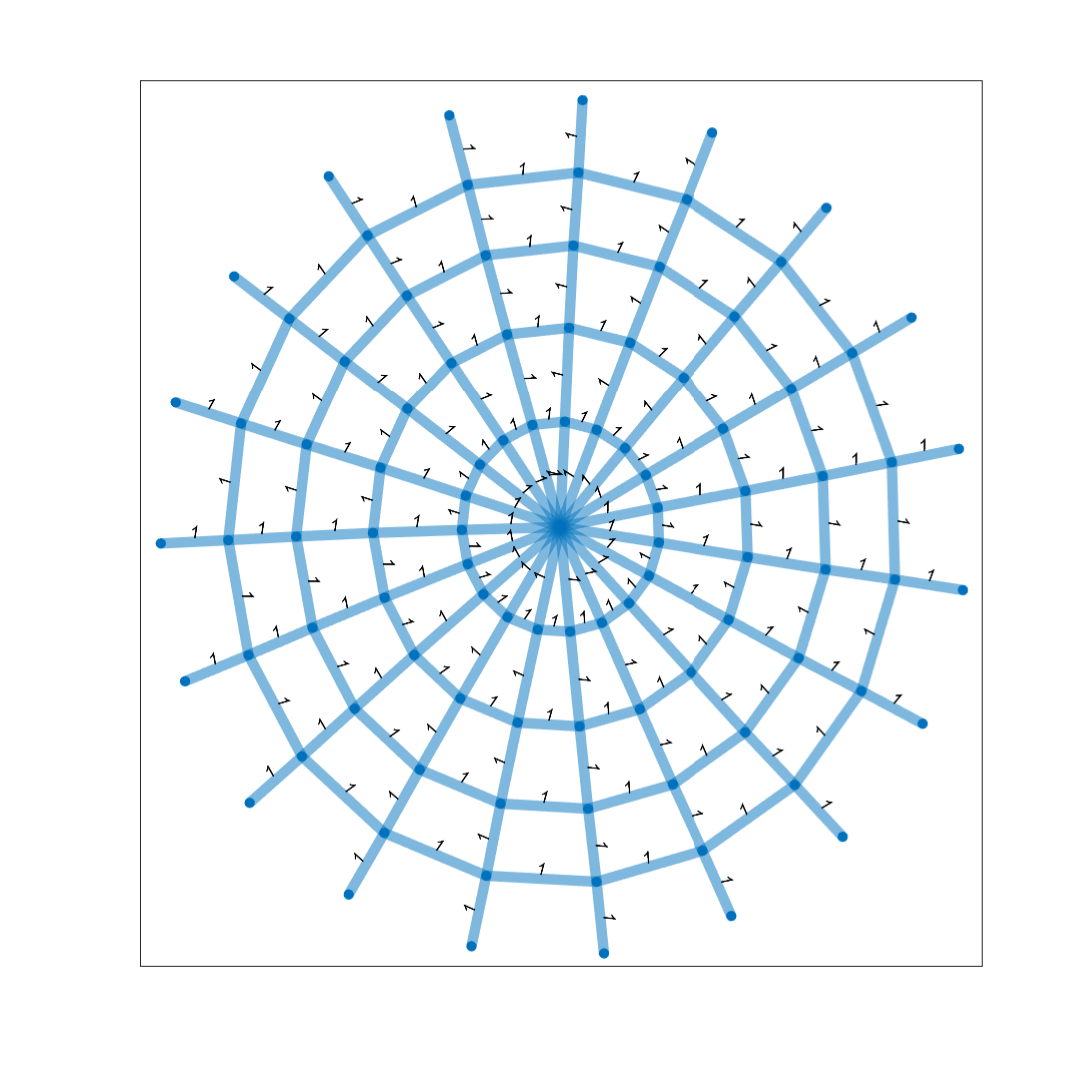}
\caption{Recovered network with $m=19$ radii in Example 1.}
\label{19}
\end{center}
\end{figure}

\begin{figure}[h!]
\begin{center}
\includegraphics[scale=0.5]{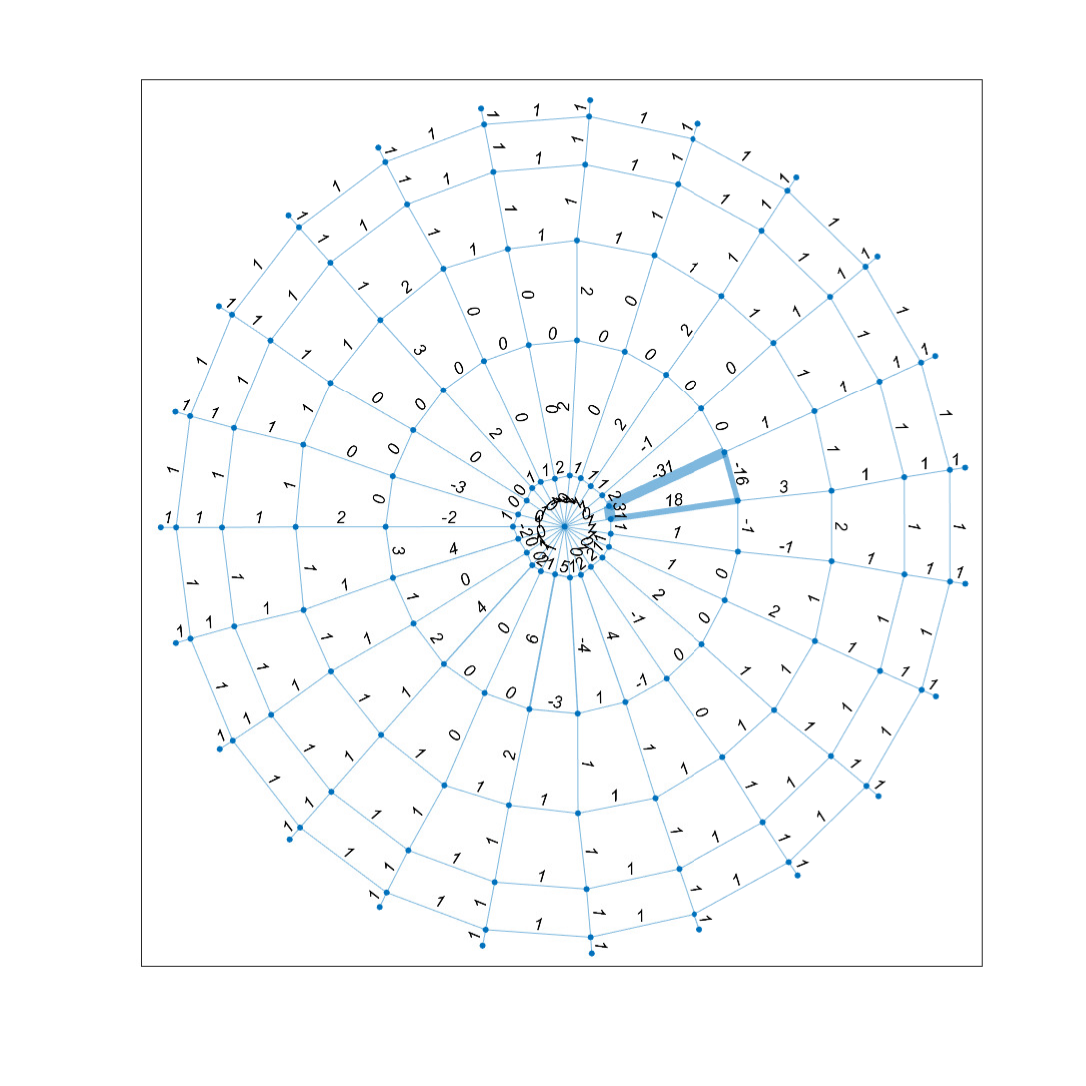}
\caption{Recovered network with $m=23$ radii
in Example 1.}
\label{23}
\end{center}
\end{figure}

We show the recovered conductance for $m=19$ and for $m=23$ in Figures \ref{19} and \ref{23}, respectively. In both figures, the width of each edge is proportional to the absolute value of the recovered conductance $c'$. For the sake of clarity, the values displayed on each edge have been rounded to the nearest integer within the graphical illustrations.

The error for $m=19$ is approximately $3.5 \cdot 10^{-3}$, and we see that the nearest integer to the value of the recovered conductance at every edge is equal to the true value $c=1$. However, the error for the next bigger network, the one with $m=23$, is approximately $5.2 \cdot 10$. We can see that the value of the recovered conductance is very far from $1$ and in some cases even negative,  especially in edges that are far from the boundary. For example there is an edge with conductance close to $-31$.

As a conclusion of the performed tests,  the recovery of the conductance of a well-connected spider network is unstable except for small networks. Moreover, the big discrepancies appear on edges that are far away from the boundary. This situation is analogous to the one that appears in the continuous Calderon's inverse conductivity problem. It is well-known that this problem is severely ill-posed, and the instabilities increase as we move farther away from the boundary.  However, in some cases, the inverse conductivity problem becomes stable if certain {\it a priori} conditions on the conductivity are known. For instance, in \cite{AV05} 
it was shown that if the conductivity is piecewise constant with a bounded number of unknown values, then the problem becomes Lipschitz stable. In the next section, we translate this hypothesis to the discrete setting.

\section{Stable recovery: discrete piecewise constant conductances hypothesis}

In the discrete case,  we define the conductance as being {\it piecewise constant} on a partition $E= E_1 \sqcup \cdots\sqcup E_s$ if it is constant on each $E_i$. Of course, as the  number of edges is finite, the conductance is inherently piecewise constant. In this work, we understand that the piecewise constant hypothesis holds if and only if $s$, the number of subsets in the partition, satisfies $s \ll |E|$. Nonetheless, when considering the extreme scenario where $s=1$ it has been observed that the conductance recovery method leads to instabilities, as demonstrated in the preceding section. Consequently, it becomes imperative to develop alternative algorithms that ensure stability. Our proposal is to formulate the inverse problem as a polynomial optimization problem 
that includes the deviation with respect  to the hypothesis of piecewise constant conductances as a penalty.

\subsection{Polynomial optimization problem}

Let $\Gamma=(\bar F,c)$ be a well-connected spider network with Laplacian matrix ${\sf L}$ and response  matrix ${\sf N}$. We denote by ${\sf e}^{t}$ the $t$-th vector of the canonical basis in $\mathbb{R}^m$, with $m=|\delta (F)|$. We define: 
 \begin{equation}\label{vf}
{\sf v}_F^{t}= {\sf L}(F;F)^{-1} \cdot  {\sf L}(\delta (F); F)^T{\sf e}^{t},
\end{equation}
which is the unique vector such that
\begin{equation}\label{eq1}
{\sf L} \begin{pmatrix}  
	{\sf e}^{t} \\
	{\sf v}_F^{t} 
\end{pmatrix}  = \begin{pmatrix}
{\sf N}(\delta(F); \{t\}) \\
	0 
\end{pmatrix};
\end{equation}
that is, ${\sf v}_F^{t}$ is the vector whose components are the potential at the nodes in $F$ when there is no current injected in $F$ and we set the potential at the $t$-th boundary vertex to $1$ and at the rest of nodes of the boundary to $0$. We denote by $v_j^t$ the $j$-th component of ${\sf v}^{t}:=\begin{pmatrix}
	{\sf e}^{t} \\
	{\sf v}_F^{t} 
\end{pmatrix} \in \mathbb{R}^n$, with $n=|\bar{F}|$.

The polynomial optimization problem that we raise to solve the discrete inverse problem can be state as follows.

\noindent {\bf Main Problem.} {\em Given a well-connected spider network $\Gamma=(\bar F,c)$ with known set of edges $E$ but unknown conductance, a response matrix ${\sf N}$, a partition $E= E_1 \sqcup \cdots\sqcup E_s$ and a penalty parameter $\mu\geq 0$; determine  values of the variables
\begin{itemize}
  
\item $c_{jk}$ for all $(j,k)\in E$;

\item $c_{j}$ for all $j=1,\ldots,s$; 

\item $v_j^t$ for all $j=m+1,\ldots, n$ and $t=1,\ldots,m,$ 
  
\end{itemize}	
\noindent which minimize the objective function
\begin{equation}\label{obj}
\hspace{-.8cm}p= \sum_{t=1}^{m}{\sum_{j=1}^{m}{\left(\sum_{k=1}^{n}{c_{jk}\left(v_j^t-v_k^t\right)} - {\sf N}\left(\{j\}; \{t\}\right)\right)^2}}   + \mu\sum_{i=1}^{s}{ \sum_{(j,k)\in E_i}{(c_{jk}-c_i)^2}}\hspace{.8cm}
\end{equation}
subject to the constraints
\begin{equation}\label{ideal}
g_j^t:= \sum_{k=1}^{n}{c_{jk}\left(v_j^t-v_k^t\right)}=0
\end{equation}
for all $j=m+1,\ldots,n$ and $t=1,\ldots,m$; and $c_{jk}\geq 0$ for all $(j,k)\in E$.
}

In our problem, we will denote by $\Gamma'=(\bar F,c')$ the spider network  whose conductance is given by the variables $c_{jk}$, so then $c'(j,k)=c_{jk}$. From (\ref{eq1}), the $(j,t)$-th entry of the response matrix ${\sf N}'$ of $\Gamma'$ is equal to $\sum_{k=1}^{n}{c_{jk}\left(v_j^t-v_k^t\right)}$ when  $c'$ and the voltage variables $v_j^t$ with $j=m+1,\ldots, n$ satisfy equations (\ref{ideal}), which are the last $n-m$ components of (\ref{eq1}) for all $t$. The physical meaning of these equations is that there is no current in the vertices in $F$ when the potential in the network is ${\sf v}^{t}$.

For each $j=1,\ldots,s$, the variable $c_{j}$ is an unknown conductance. We denote by ${\bar c}$ the piecewise constant conductance on the partition $E= E_1 \sqcup\cdots\sqcup E_s$ which is equal to $c_{j}$ on each $E_j$.
	
The objective function (\ref{obj}) can be rewritten as 
\begin{equation}\label{obj2}
	p= || {\sf N}' - {\sf N} ||^2    + \mu || c'-\bar c ||^2,
\end{equation}
so  a solution of the main problem  minimizes the squared Frobenius norm of the difference between the response matrix ${\sf N}'$ of the recovered network and ${\sf N}$ plus a penalty term which is the squared Euclidian norm of the difference between the recovered conductance and any piecewise constant conductance on $E= E_1 \sqcup \cdots\sqcup E_s$ multiplied by the penalty parameter $\mu$. In the context of Tikhonov-like  regularisation methods the parameter $\mu$ is often called regularisation parameter (see \cite{EHN96,LMP03,R16}).  

The case $\mu \rightarrow \infty$ corresponds with enforcing the hypothesis that the recovered conductance is piecewise constant on $E= E_1 \sqcup\cdots \sqcup E_s$ and minimizing the difference between ${\sf N}'$ and ${\sf N}$. The case $\mu=0$ corresponds with minimizing the difference between ${\sf N}'$ and ${\sf N}$ ignoring the piecewise constant hypothesis.

A solution to the main problem must be a stationary point of the Lagrangian function of the optimization problem. In particular, for each $i=1,\ldots,s$, the partial derivative of $\sum_{(j,k)\in E_i}{(c_{jk}-c_i)^2}$ with respect to the variable $c_i$ must be equal to zero, so

\begin{equation}\label{ci}
c_i= \frac{1}{|E_i|}\sum_{(j,k)\in E_i}{c_{jk}}.
\end{equation}
Therefore, the $c_i$ variable; that is, the value of the  conductance on the partition can be removed from the optimization problem by  substituting in (\ref{obj}) the variables $c_i$ with (\ref{ci}).

In a well-connected spider network it is easy to check that there are $|E|=\displaystyle{m \choose 2}$ edges and $n=\displaystyle\frac{m^2+m+4}{4}$ vertices. So, the number of variables in our problem  is equal to
$$
r=|E|+m(n-m)=\frac{m(m^2-m+2)}{4}.
$$

We denote by $J=\big\langle g_j^t;  j=m+1,\ldots,n, \, \,   t=1,\ldots,m\big\rangle$ the ideal generated by the quadratic polynomials in (\ref{ideal}), and we denote by $V(J) \subset \mathbb{R}^{r}$ the real vanishing set of $J$. Then, the main problem  can be stated as finding a minimum of the quartic $p$ in $\left(\mathbb{R}_{\geq 0}^{|E|} \times \mathbb{R}^{r-|E|}\right) \cap V(J)$.

\subsection{Problem resolution}

For $\mu=0$, the main problem  has a unique solution $\hat y\in \left(\mathbb{R}_{\geq 0}^{|E|} \times \mathbb{R}^{r-|E|}\right) \cap V(J)$ with $p(\hat y)=0$, which is $c'=c$ and the value of the voltage variables is given by (\ref{vf}). For any $\mu>0$, if $c$ is piecewise constant in the partition $E= E_1 \sqcup \cdots\sqcup E_s$, then the main problem also has the same unique solution with $p(\hat y)=0$.

In the general case, with a numerical optimization method, such as an interior point algorithm, we can obtain an approximation to a minimum $y^*\in \left(\mathbb{R}_{\geq 0}^{|E|} \times \mathbb{R}^{r-|E|}\right) \cap V(J)$ of the objective function $p$. In Section \ref{res}, we will discuss several examples of how this minimum approximates the solution of our problem when the conductance is piecewise constant.

As an initial guess for the interior point algorithm, for every conductance variable, we choose the conductance of the spider with constant conductance $c^0$  whose response matrix is the closest to ${\sf N}$ in the Frobenius norm.

\begin{lemma}\label{ini}
	Given a network with boundary $\Gamma=(\bar F,c)$ and its Dirichlet-to-Neumann matrix ${\sf N}$, then, the network $\Gamma^0=(\bar F,c^0)$ with the same set of edges as $\Gamma$ and constant conductance $c^0\ge 0$ at all edges whose Dirichlet-to-Neumann matrix is the closest to ${\sf N}$ in the Frobenius norm satisfies that $c^0$ is the solution to the following nonnegative least squares (NNLS) problem
$$
\min\limits_{c^0\ge 0}\left\{||{\sf N} - c^0 ({\sf I} - {\sf L}'(\delta (F); F) {\sf L}'(F; F)^{-1} {\sf L}'(\delta (F); F)^T)||^2\right\},
$$
where ${\sf L}'$ is the (unweighted) combinatorial Laplacian of the spider graph.
\end{lemma}

NNLS problems are convex \cite{SH13}, so every critical point is a global minimum. We propose here to solve the above problem also with an interior point method, which has been used successfully to solve NNLS problems in the context of inverse problems in electrical networks (see \cite{S23}).

\begin{remark}\label{inir}
In the case of a well-connected spider, if $m\geq 7$, the NNLS problem in Lemma \ref{ini} becomes
$$
\min\limits_{c^0\ge 0}\left\{||{\sf N} - c^0 ({\sf I} - {\sf G}(\delta(F^c); \delta(F^c))||^2\right\}
$$
where ${\sf G}={\sf L}'(F; F)^{-1},$ 
 and if $m=3$, the NNLS problem becomes
$$
\min\limits_{c^0\ge 0}\left\Vert {\sf N} - \frac{c^0}{3}\begin{pmatrix}
	2 & -1 & -1 \\
	-1 & 2 & -1 \\
	-1 & -1 & 2 
\end{pmatrix} \right\Vert^2.
$$
\end{remark}

\begin{remark}
As an initial guess for every voltage variable $v_j^t$, we choose the voltage at node $j$, for the spider $\Gamma^0$  with constant conductance,  which gives zero current in $F$ when the voltage in $\delta(F)$ is equal to ${\sf e}^t$; that is, for every $t=1,\ldots,m$, the initial guesses of the voltage variables  are given by
\begin{equation}\label{vf0}
({\sf v}_F^{t})^0={\sf L}'(F; F)^{-1} \cdot  {\sf L}'(\delta (F); F)^T{\sf e}^{t}.
\end{equation}
\end{remark}

\section{Numerical results} \label{res}
In this section, we present numerical examples for the recovery of piecewise constant conductances on  well-connected spider networks, using the method explained in the previous section. Our approach effectively avoids the typical instabilities associated with these types of problems.

We write the optimization problem raised in the previous section in MATLAB using Casadi \cite{Ca23}, an open-source software tool that provides a symbolic framework suited for numerical optimization, and we obtain an approximate solution of it using the interior-point solver IPOPT \cite{Ip23}, an
open-source software package for large-scale nonlinear optimization. Moreover, to obtain the initial guess of the variables we also use IPOPT. The tolerance used in IPOPT is equal to $10^{-8}$. 
\begin{example}\label{estable}
For each $m=7, 11, 15, 19, 23, 27, 31$ and $35$, and for each $s=1,\ldots, 10$, we generate $10$ well-connected spider networks with $m$ radii, each with a conductance which is piecewise constant in a (possibly different) random partition with $s$ subsets. For the case $m=3$, we do the same, but only for $s=1, 2$ and $3$, because the total number of edges of the graph is $3$.

For each combination of values of $m$ and $s$, for each of the $10$ networks, the value of the conductance at each subset of the partition is sampled from the uniform distribution in the interval $[1,100]$. The total number of tested networks is $830$.

\end{example}

For each network with conductance $c$, we compute its response matrix ${\sf N}$, we recover the conductance $c'$, we compute the error $||c'-c||$, and we compute the quotient $\displaystyle\frac{||c'-c||}{||{\sf N}'-{\sf N}||}$, which can be seen as an approximation of the Lipschitz constant for the case of piecewise constant conductance. In Table \ref{tabla_c}, we show the maximum error in the recovered conductance in the $10$ networks with each combination of values of $m$ and $s$. Analogously, we show in Table \ref{tabla_k} the maximum $\displaystyle\frac{||c'-c||}{||{\sf N}'-{\sf N}||}$ in the $10$ networks with each combination of values of $m$ and $s$.

\begin{table}[h!]
{\small{
\caption{\label{tabla_c}Maximum error in the recovered conductance with one significative digit.}
\begin{tabular}{@{}llllllllllll}
\\[-3ex]
\noalign{\hrule height 1pt} 
$m\backslash s$ & 1 & 2  & 3 & 4  & 5 & 6 & 7 & 8  & 9 & 10\\
\hline
3 & $2\cdot 10^{-10}$ &$8\cdot 10^{-10}$& $1\cdot 10^{-9}$  &   & &  &    &&   & \\
7 & $1\cdot 10^{-8}$ &$1\cdot 10^{-10}$& $6\cdot 10^{-6}$  & $1\cdot 10^{-9}$  &$2\cdot 10^{-8}$& $1\cdot 10^{-8}$ & $4\cdot 10^{-7}$   &$9\cdot 10^{-8}$& $4\cdot 10^{-7}$  & $1\cdot 10^{-6}$\\
11 & $3\cdot 10^{-9}$ &$5\cdot 10^{-9}$& $4\cdot 10^{-10}$  & $3\cdot 10^{-8}$  &$2\cdot 10^{-8}$& $9\cdot 10^{-9}$ & $9\cdot 10^{-8}$   &$3\cdot 10^{-8}$& $3\cdot 10^{-7}$  & $4\cdot 10^{-7}$\\
15 & $5\cdot 10^{-8}$ &$4\cdot 10^{-10}$& $3\cdot 10^{-8}$  & $2\cdot 10^{-9}$  &$8\cdot 10^{-8}$& $6\cdot 10^{-9}$ & $2\cdot 10^{-7}$   &$2\cdot 10^{-8}$& $1\cdot 10^{-7}$  & $8\cdot 10^{-8}$\\
19 & $2\cdot 10^{-8}$ &$2\cdot 10^{-9}$& $2\cdot 10^{-8}$  & $2\cdot 10^{-9}$  &$2\cdot 10^{-8}$& $2\cdot 10^{-8}$ & $2\cdot 10^{-8}$   &$4\cdot 10^{-8}$& $2\cdot 10^{-7}$  & $1\cdot 10^{-6}$\\
23 & $3\cdot 10^{-8}$ &$1\cdot 10^{-8}$& $5\cdot 10^{-9}$  & $2\cdot 10^{-8}$  &$1\cdot 10^{-8}$& $7\cdot 10^{-9}$ & $1\cdot 10^{-6}$   &$3\cdot 10^{-7}$& $1\cdot 10^{-7}$  & $9\cdot 10^{-8}$\\
27 & $2\cdot 10^{-7}$ &$4\cdot 10^{-8}$& $9\cdot 10^{-8}$  & $1\cdot 10^{-6}$  &$2\cdot 10^{-7}$& $6\cdot 10^{-8}$ & $5\cdot 10^{-7}$   &$9\cdot 10^{-8}$& $1\cdot 10^{-7}$  & $4\cdot 10^{-7}$\\
31 & $2\cdot 10^{-7}$ &$2\cdot 10^{-8}$& $4\cdot 10^{-8}$  & $3\cdot 10^{-7}$  &$7\cdot 10^{-9}$& $3\cdot 10^{-8}$ & $3\cdot 10^{-8}$   &$1\cdot 10^{-7}$& $4\cdot 10^{-8}$  & $2\cdot 10^{-7}$\\
35 & $1\cdot 10^{-7}$ &$3\cdot 10^{-8}$& $1\cdot 10^{-7}$  & $7\cdot 10^{-8}$  &$1\cdot 10^{-7}$& $7\cdot 10^{-8}$ & $7\cdot 10^{-7}$   &$1\cdot 10^{-7}$& $1\cdot 10^{-7}$  & $8\cdot 10^{-7}$\\
\noalign{\hrule height 1pt} 
\end{tabular}}}
\end{table}
\normalsize

\begin{table}[h!]
{{
\caption{\label{tabla_k}Maximum $\frac{||c'-c||}{||{\sf N}'-{\sf N}||}$ in the recovered conductance.}
\begin{tabular}{@{}llllllllllll}
\\[-3ex]
\noalign{\hrule height 1pt} 
$m\backslash s$ & 1 & 2  & 3 & 4  & 5 & 6 & 7 & 8  & 9 & 10\\ 
\hline
3 & $1\cdot 10^{0}$ &$4\cdot 10^{-1}$& $7\cdot 10^{0}$  &   & &  &    &&   & \\
7 & $3\cdot 10^{0}$ &$4\cdot 10^{0}$& $2\cdot 10^{2}$  & $3\cdot 10^{1}$  &$1\cdot 10^{2}$& $2\cdot 10^{2}$ & $8\cdot 10^{1}$   &$6\cdot 10^{2}$& $1\cdot 10^{2}$  & $7\cdot 10^{2}$\\
11 & $4\cdot 10^{0}$ &$8\cdot 10^{0}$& $9\cdot 10^{0}$  & $1\cdot 10^{1}$  &$7\cdot 10^{1}$& $7\cdot 10^{1}$ & $3\cdot 10^{1}$   &$7\cdot 10^{1}$& $8\cdot 10^{2}$  & $1\cdot 10^{2}$\\
15 & $5\cdot 10^{0}$ &$5\cdot 10^{0}$& $4\cdot 10^{1}$  & $1\cdot 10^{1}$  &$5\cdot 10^{1}$& $1\cdot 10^{1}$ & $5\cdot 10^{1}$   &$6\cdot 10^{1}$& $1\cdot 10^{2}$  & $8\cdot 10^{1}$\\
19 & $5\cdot 10^{0}$ &$6\cdot 10^{0}$& $7\cdot 10^{0}$  & $1\cdot 10^{1}$  &$8\cdot 10^{0}$& $3\cdot 10^{1}$ & $6\cdot 10^{1}$   &$9\cdot 10^{1}$& $2\cdot 10^{2}$  & $8\cdot 10^{1}$\\
23 & $6\cdot 10^{0}$ &$7\cdot 10^{0}$& $1\cdot 10^{1}$  & $4\cdot 10^{1}$  &$1\cdot 10^{1}$& $2\cdot 10^{1}$ & $2\cdot 10^{1}$   &$3\cdot 10^{1}$& $2\cdot 10^{1}$  & $2\cdot 10^{1}$\\
27 & $7\cdot 10^{0}$ &$7\cdot 10^{0}$& $3\cdot 10^{1}$  & $8\cdot 10^{0}$  &$1\cdot 10^{1}$& $2\cdot 10^{1}$ & $8\cdot 10^{1}$   &$2\cdot 10^{1}$& $8\cdot 10^{1}$  & $7\cdot 10^{1}$\\
31 & $7\cdot 10^{0}$ &$8\cdot 10^{0}$& $9\cdot 10^{0}$  & $9\cdot 10^{0}$  &$2\cdot 10^{1}$& $2\cdot 10^{1}$ & $6\cdot 10^{1}$   &$5\cdot 10^{1}$& $6\cdot 10^{1}$  & $4\cdot 10^{1}$\\
35 & $8\cdot 10^{0}$ &$9\cdot 10^{0}$& $1\cdot 10^{1}$  & $1\cdot 10^{1}$  &$1\cdot 10^{1}$& $1\cdot 10^{1}$ & $1\cdot 10^{1}$   &$1\cdot 10^{1}$& $3\cdot 10^{1}$  & $1\cdot 10^{1}$\\
\noalign{\hrule height 1pt} 
\end{tabular}}}
\end{table}
\normalsize

To recover the conductance of all networks, we choose a value of $\mu=1$ as penalty parameter in our problem. With this value, in all the experiments, we obtain a solution such that the evaluation of the quartic $p$ at it is lower than $1.3746\cdot 10^{-13}$, the evaluation of the deviation with respect to the hypothesis of being piecewise constant, $||c'-\bar c||^2$ is lower than $9.1998\cdot 10^{-14}$, and the deviation in the response matrix with respect to the data, $||{\sf N}'-{\sf N}||^2$, is lower than $4.5463\cdot 10^{-14}$.

We see in Table \ref{tabla_c} that we recover the conductance of all networks with an error lower than $6.0250\cdot 10^{-6}$, which is much lower than the norm of the real conductances. In Table \ref{tabla_k}, we see that in all cases the error in the conductance is lower than $833$ times the error in the response matrix with respect to the data, $||{\sf N}'-{\sf N}||$. 

As a note, we have done the same experiment for different positive values of $\mu$, recovering in all cases the conductances with very low error. Nevertheless, as $\mu \rightarrow 0$, the recovery becomes unstable for large networks, because we obtain solutions such that the evaluation of the quartic $p= ||{\sf N}' - {\sf N} ||^2    + \mu || c'-\bar c ||^2$ is almost zero, but the value $|| c'-\bar c ||^2$ is very high, so the conductance is far from being piecewise constant. 

In all $830$ cases, the fact that the deviation with respect to the piecewise constant hypothesis is very low gives us a stable recovery of the conductances. In Calderon's inverse problem the Lipschitz stability constant 
depends on the partition in which the conductivity is piecewise constant, as well as in the bounds of the conductivity and the set in which the conductivity is defined. It is expected that the constant will diverge as the number of subsets in the partition goes to infinity \cite{AV05}. 

By analogy, we expect that the stability of the problem depends on $m$, $s$, the partition in which the conductance is piecewise constant and the bounds of the conductance. We also expect the problem becomes unstable as $s \rightarrow |E|$, because in the limit $s=|E|$, each edge is in a different subset of the partition, so the piecewise constant conductance hypothesis is meaningless. As the partition is different in general for the $10$ networks with the same value of $m$ and $s$, we get quite different errors recovering them, but we see that in general, the maximum of the ratio $\displaystyle\frac{||c'-c||}{||{\sf N}'-{\sf N}||}$ increases with $s$, as expected. For example, in the case $m=7$ and $s=10$, $s$ is almost half of the total number of edges ($|E|=21$), and the ratio is equal to $7\cdot 10^2$. However, for $m=35$ and $s=10,$ that represents almost the $2\%$ of the edges,  the ratio equals $10.$

Lastly, we will show an example of recovery of piecewise constant conductances in networks with the same value of radii $m=11$ and all possible number of subsets of the partition $s$ to see how the stability varies with $s$.

\begin{example}\label{ej11}
We set $m=11$, and for each $s=1,\ldots, |E|=55$, we generate $100$ different spider networks with $11$ radii and piecewise constant conductance whose  values  at each subset of the partition is sampled from the uniform distribution in the interval $[1,100]$. 
\end{example}

\begin{figure}[ht!]
\begin{center}
\includegraphics[scale=0.6]{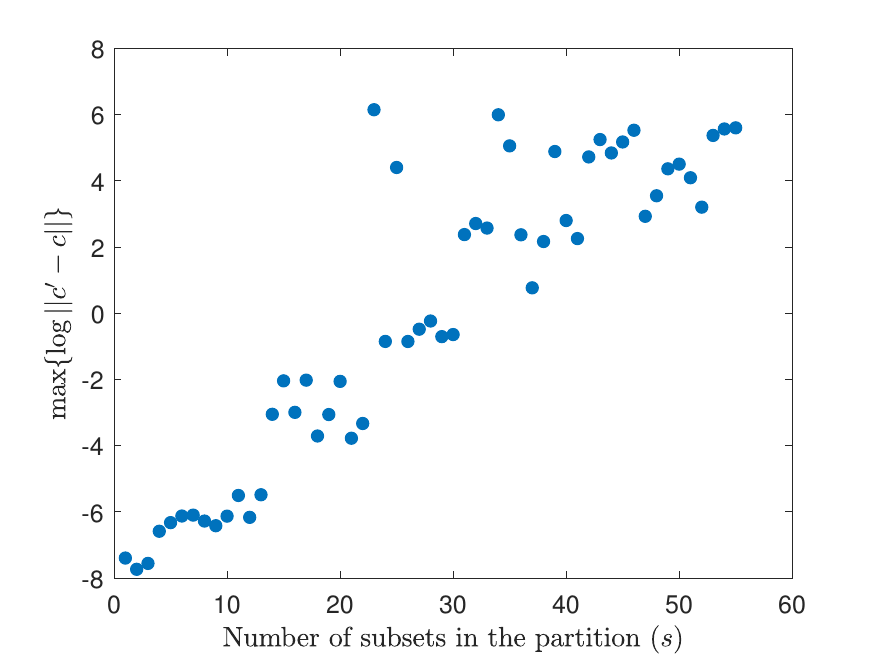}
\caption{Maximum of the logarithm of the error in the recovered conductances of $100$ networks with $m=11$ and $s=1,\ldots,55$. (Example \ref{ej11}).}
\label{c11}
\end{center}
\end{figure}

\begin{figure}[ht!]
\begin{center}
\includegraphics[scale=.65]{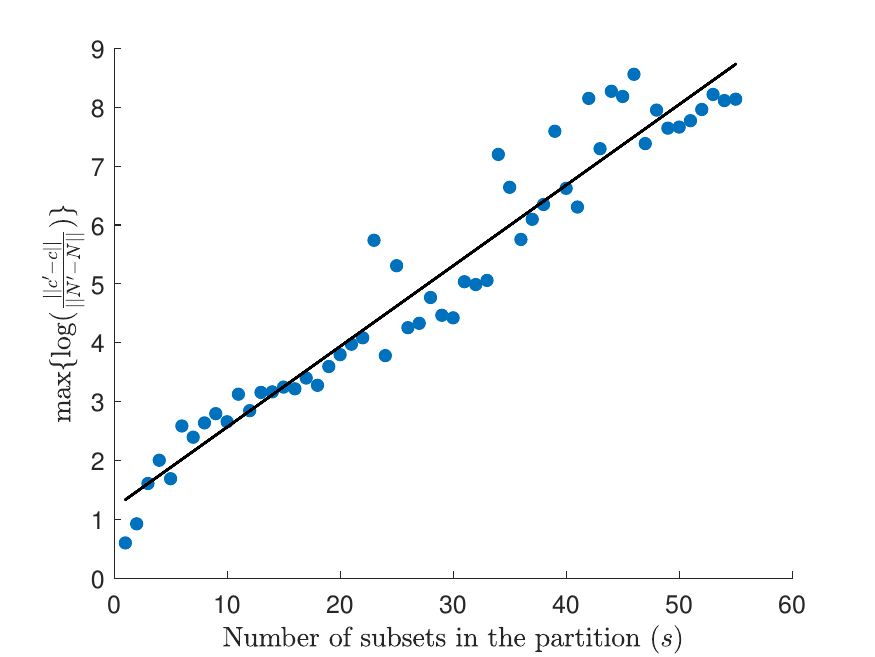}
\caption{Maximum  $\log\left(\frac{||c'-c||}{||{\sf N}'-{\sf N}||}\right)$ in  recovered conductance of $100$ networks with $m=11$ and $s=1,\ldots,55$. (Example \ref{ej11}).}
\label{k11}
\end{center}
\end{figure}

As in Example \ref{estable}, we set $\mu=1$, we compute its response matrix ${\sf N}$, we recover the conductance $c'$, we compute the error $||c'-c||$, and we compute the quotient $\displaystyle\frac{||c'-c||}{||{\sf N}'-{\sf N}||}$. In Figure \ref{c11}, we show the maximum error in the recovered conductance in the $100$ networks with each value of $s$. Analogously, we show in Figure \ref{k11}, the maximum $\displaystyle\frac{||c'-c||}{||{\sf N}'-{\sf N}||}$ in the $100$ networks with each value of $s$.

In this example, we see that there is a strong linear dependence between the maximum of the ratio $\displaystyle\frac{||c'-c||}{||{\sf N}'-{\sf N}||}$ and $s$, which is in agreement with the results obtained in \cite{AV05,R06} about  exponentially behaviour  of the Lipschitz constant with respect to the number of regions where the conductivity is constant. The line in Figure \ref{k11} shows the linear regression between the number of subsets in a partition and the maximum ratio. It has been computed using Matlab and the results show a value of R-squared equal to $0.938$ with a p-value of $1.22\cdot 10^{-33}.$

In every case, we are able to recover a network with a low value of $||{\sf N}'-{\sf N}||$, so we see that the error $||c'-c||$ is low when $s \ll |E|$. In particular, we see that for values of $s$ between $1$ and $13$, the error in the recovered conductances is very low in all the networks. Besides,  for values of $s$ between $14$ and  $22$, the error is greater in some cases, but smaller than the norm of the real conductances. Finally,   for $s>23$, there are networks in which the recovered conductances have a great error.

\printcredits

\bibliographystyle{cas-model2-names}

\bibliography{biblio}
\end{document}